\documentclass[a4paper,11pt]{article}
\usepackage[latin1]{inputenc}
\usepackage{epsfig}
\usepackage{amsmath,amssymb}
\usepackage{amsfonts}
\usepackage{indentfirst}

\title{Some statistics on permutations avoiding generalized patterns\footnote{This work was partially
supported by MIUR project: \emph{Automi e linguaggi formali:
aspetti matematici e applicativi}.}}
\date{}

\author{Antonio Bernini\thanks{Università di Firenze, Dipartimento di Sistemi e Informatica,
viale Morgagni 65, 50134 Firenze, Italy \tt{bernini@dsi.unifi.it,
ferrari@math.unifi.it}} \and Mathilde Bouvel\thanks{\'Ecole
Normale Supérieure de Cachan, Département informatique, 61 av. du
Président Wilson, 94235, Cachan cedèx, France
\tt{mathilde.bouvel@dptinfo.ens-cachan.fr}}
 \and Luca
Ferrari$^{\dag}$}

\begin{document}
\hyphenation{ge-ne-ra-li-zed re-wri-ting}

\maketitle

\begin{abstract}
In the last decade a huge amount of articles has been published
studying pattern avoidance on permutations. From the point of view
of enumeration, typically one tries to count permutations avoiding
certain patterns according to their lengths. Here we tackle the
problem of refining this enumeration by considering the statistics
``first/last entry". We give complete results for every
generalized patterns of type $(1,2)$ or $(2,1)$ as well as for
some cases of permutations avoiding a pair  of generalized
patterns of the above types.
\end{abstract}

\section{Introduction}
Let $\pi\in S_n$ and $\tau\in S_k$, where $S_t$ denotes the
symmetric group on $[t]=\{1,2,\ldots,t\}$. We say that $\pi$
\emph{avoids} $\tau$ if there are no subsequences
$\pi_{i_1}\pi_{i_2}\ldots \pi_{i_k}$ with $1\leq
i_1<i_2<\ldots<i_k\leq n$ which are order-isomorphic to $\tau$,
that is having the entries $\pi_{i_j}$ in the same relative order
of the entries of $\tau$. The permutation $\tau $ is called a
(classical) \emph{pattern}. We denote the set of all
$\tau$-avoiding permutations of $S_n$ with $S_n(\tau)$. In
\cite{bab-stein}, \emph{generalized patterns} were introduced to
study the Mahonian statistics on permutations. They are obtained
by inserting one or more dashes among the entries of $\tau$. A
pattern $\tau=\tau_1-\tau_2-\ldots -\tau_k$ with $k-1$ dashes is
called of type $(|\tau_1|,|\tau_2|,\ldots |\tau_k|)$, where
$|\tau|$ is the length of $\tau$. For instance, $\tau=13-26-574$
is a pattern of type $(2,2,3)$. A classical pattern of length $k$
can be seen as a pattern of type $\underbrace{(1,1,\ldots
,1)}_{k}$, assuming that a dash is inserted, but not showed,
between each pair of consecutive elements of the classical
pattern. If $\tau\in S_3$, then generalized patterns deriving from
$\tau$ are of type $(1,2)$ or $(2,1)$ according to the number of
elements preceding and following the dash and they are collected
in the set
$$
            \begin{array}{ll}
            \mathcal{P}  =& \{1-23,12-3,1-32,13-2,3-12,31-2,2-13,21-3,\\
                  &\mbox{  } \mbox{  } 2-31,23-1,3-21,32-1\}.\\
            \end{array}
$$
A permutation $\pi$ contains a generalized pattern $p\in
\mathcal{P}$ if adjacent elements in $p$ are also adjacent in
$\pi$. For example $\pi=7256134$ contains the generalized pattern
$13-2$ in its subsequence $\pi_2\pi_3\pi_6=253$. Note that it does
not contain the pattern $1-32$, but it contains the classical
pattern $132$ in the subsequences $\pi_2\pi_4\pi_6$ and
$\pi_2\pi_4\pi_7$.

Let $\pi$ be a permutation of $S_n$, then the \emph{reverse} and
\emph{complement} permutations $\pi^r$ and $\pi^c$, respectively,
are defined as follows: $\pi^r_i=\pi_{n+1-i}$ and
$\pi^c_i=n+1-\pi_i$, for $i=1,\ldots ,n$. We can define the
reverse and the complement also for a pattern, regarding the dash
as a particular entry in reversing $\pi$ and leaving it in the
same position when the complement $\pi^c$ is performed. So, for
example, if $p=1-32$, then $p^r=23-1$ and $p^c=3-12$. Considering
the composition of the reverse and the complement, it is easily
seen that $p^{cr}=p^{rc}$. The set $\{p,p^r, p^c, p^{cr}\}$ is
called the \emph{symmetry class} of $p$. Observe that
$|S_n(p)|=|S_n(p^r)|=|S_n(p^c)|=|S_n(p^{rc})|$.

The twelve generalized patterns of $\mathcal P$ are organized in
three symmetry classes : $\{1-23,32-1,3-21,12-3\}$, $\{3-12,21-3,
1-32,23-1\}$ and $\{2-13,31-2,2-31,13-2\}$. If $p$ and $p'$ are
two patterns such that $|S_n(p)|=|S_n(p')|$, then $p$ and $p'$ are
said to be in the same \emph{Wilf class} \cite{man}. Since in
\cite{clae} it is shown that
\begin{itemize}
    \item $|S_n(p)|=B_n$,\\ for $p\in\{1-23,32-1,3-21,12-3\}\bigcup
    \{3-12,21-3,1-32,23-1\}$
    \item $|S_n(p)|=C_n$,\\ for $p\in\{2-13,31-2,2-31,13-2\},$
\end{itemize}
where $B_n$ and $C_n$ are the $n$-th Bell and Catalan numbers,
respectively, then we can say that $\mathcal{P}$ is organized in
two Wilf classes: $\{1-23,32-1,3-21,12-3,3-12,21-3,1-32,23-1\}$
and $\{2-13,31-2,2-31,13-2\}.$

In this work, we refine some enumerative results on $S(p)$, $p\in
\mathcal{P}$, namely we count $p$-avoiding permutations, for each
$p$, according to their length and the value of their first or
last entry. Next we solve the same problem for some classes of
permutations of the kind $S(p,q)$, $p,q \in \mathcal{P}$, and we
conclude by proposing to tackle this problem for any remaining
pair of generalized patterns of $\mathcal{P}$.

Our results are achieved by using the ECO method together with a
graphical representation of permutations. In the following we only
briefly recall the ECO construction for (patterns avoiding)
permutations, for more details we refer the reader to
\cite{bar-del-per-pin} and \cite{ber-fer-pin}.

\begin{figure}
\begin{center}
\includegraphics[scale=0.4]{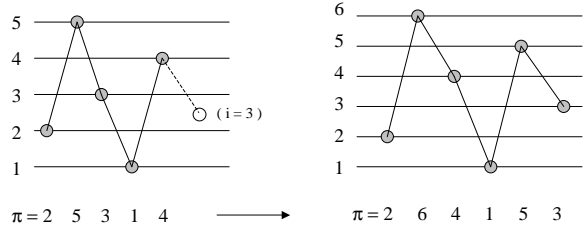}
\end{center}
\caption{an ECO construction for permutations}\label{perm}
\end{figure}

Any permutation of length $n$ can be visualized using a path-like
representation, as in Figure \ref{perm}. Note that the plane is
divided in $n+1$ strips by the $n$ horizontal lines which are
numerated from $1$ to $n$, starting from bottom (in the sequel, we
refer to these strips as ``regions": region $i$ is included
between line $i-1$ and line $i$, whereas region 1 is the one below
line 1 and region $n+1$ is the one above line $n$). Each entry of
the permutation is represented as a ``node" lying on the line
corresponding to its value, . If $\pi\in S_n$, then $n+1$
permutations belonging to $S_{n+1}$ can be obtained by inserting a
new node in each region of the plane. If we wish to generate the
permutations in $S_{n+1}(P)$ obtained in such a way from $\pi\in
S_n(P)$, where $P$ is a set of forbidden patterns, then the
regions the last node can be inserted in form a subset of all the
$n+1$ possible regions; in the framework of the ECO method they
are called \emph{active sites} \cite{bar-del-per-pin}. A
remarkable feature of this construction is that, if $\pi\in
S_n(P)$, then $\pi'\in S_{n+1}$ (which is obtained from $\pi$ by
inserting the last node in one of the regions) does not contain
the patterns specified in $P$ in its entries $\pi'_j$ with
$j=1,\ldots,n$, otherwise $\pi$ itself would contain some pattern
of $P$. So, to decide if a region $i$ is an active site or not, we
just have to check those generalized patterns the last node is
involved in.

\section{The symmetry class $\{1-23, 32-1, 3-21, 12-3 \}$}

\subsection{ECO construction and generating tree of
$S(1-23)$}\label{eco.tree.1-23}

\begin{figure}
\begin{center}
\includegraphics[scale=0.4]{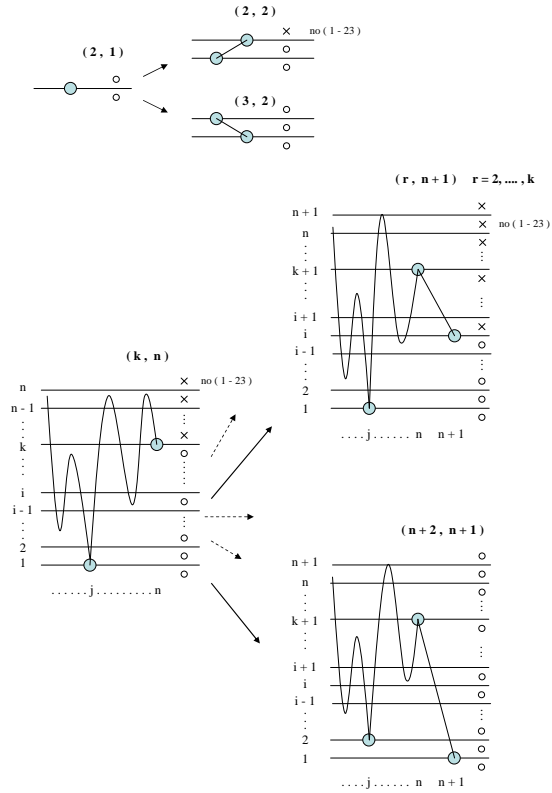}
\end{center}
\caption{ECO construction of $S(1-23)$}\label{fig:ECO1-23}
\end{figure}

Let $\pi \in S_n(1-23)$. If $\pi_n=k\neq 1$, then $\pi$ generates
$k$ permutations $\pi^{(i)}\in S_{n+1}(1-23)$, $i=1,2,\ldots ,k$,
by inserting a new node in region $i$. If $\pi_n=1$, then $\pi$
generates $n+1$ permutations by inserting a new node in any
region. Note that in this case the number of sons of $\pi$ is
determined by the length of $\pi$. If $\pi^{(r)}\in S_{n+1}(1-23)$
denotes the permutation of $S_{n+1}(1-23)$ derived from $\pi \in
S_n{(1-23)}$ by inserting the last node in region $r$, it is
easily seen that $\pi^{(1)}$ generates, in turns, $n+2$
permutations, whereas $\pi^{(r)}$, $r\neq 1$, produces $r$
permutations of $S_{n+2}(1-23)$. This ECO construction can be
represented as in Figure \ref{fig:ECO1-23} and, if we label with
$(k,n)$ each permutation of $S_n(1-23)$ having $k$ active sites,
it can be encoded by the following succession rule:
$$\Omega :
\left\{
\begin{array}{l}
  (2,1) \\
  (k,n) \rightsquigarrow (2,n+1) (3,n+1) \cdots (k,n+1)
(n+2,n+1) \quad .\\
\end{array}
\right.
$$

\begin{figure}
\begin{center}
\includegraphics[scale=0.4]{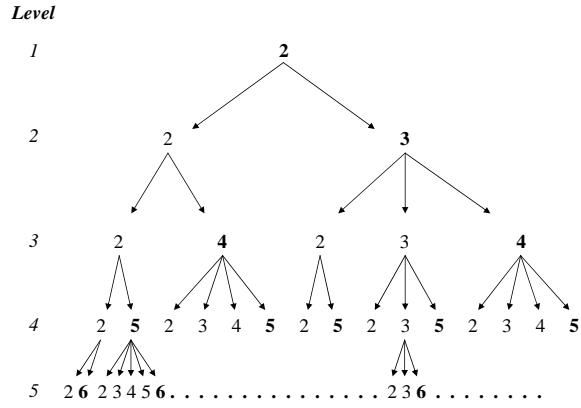}
\end{center}
\caption{the generating tree of $S(1-23)$}\label{fig:gen_tree}
\end{figure}

We now wish to draw the generating tree related the the previous
succession rule. For the sake of simplicity and for reasons that
will become clear later, we choose to label the nodes of the
generating tree using the number of their sons, which correspond
to the first element of each label of the succession rule. In
Figure \ref{fig:gen_tree} we have depicted the first levels of the
generating tree of $S(1-23)$. Here the labels in bold character
correspond to the labels of the kind $(n+1,n)$ in the succession
rule. Observe that the production of each label depends on its
level in the generating tree.

\subsection{Distribution according to the length and the last value}
\label{distr1-23}

Starting from the generating tree of Figure \ref{fig:gen_tree}, we
can consider the matrix $M=(m_{ij})_{i,j\geq 1}$ where $m_{i,j}$
is the number of labels $j+1$ at level $i$ in the generating tree:

\scriptsize
$$
M=
\left(%
\begin{array}{ccccccc}
  1 & 0 & 0 & 0 & 0 & 0 & \vdots \\
  1 & 1 & 0 & 0 & 0 & 0 & \vdots \\
  2 & 1 & 2 & 0 & 0 & 0 & \vdots \\
  5 & 3 & 2 & 5 & 0 & 0 & \vdots \\
  15 & 10 & 7 & 5 & 15 & 0 & \vdots \\
  52 & 37 & 27 & 20 & 15 & 52 & \vdots \\
  \cdots & \cdots & \cdots & \cdots & \cdots & \cdots & \ddots \\
\end{array}%
\right)
$$

\normalsize The above matrix $M$ is called the \emph{ECO matrix}
of the rule $\Omega$, according to \cite{deu-fer-rin}. It is
easily seen that $M$ can be recursively described as follows:
\begin{enumerate}
    \item $m_{1,1}=1$ (the minimal permutation $\pi=1$ has two sons);
    \item $m_{n,k}=0$ if $k>n$ (each permutation of length $n$ has at most $n$ sons);
    \item $m_{n,k}=\sum_{i=k}^{n-1}{m_{n-1,i}}$ if $k<n$ (this derives directly from
    the recursive interpretation of the previous succession rule);
    \item $m_{n,n}=m_{n,1}$ (each permutation of length
    $n-1$ produces precisely one son having label 2 and precisely one son having label $n+1$).
\end{enumerate}

Since $m_{n,1}(=m_{n,n})$ is the sum of all the elements in the
$(n-1)$-th row (for $n>1$), this entry records the total number of
$(1-23)$-avoiding permutations of length $n-1$. In other words,
$m_{n,1}=B_{n-1}$.

\bigskip
Moreover, from a careful inspection of $M$, we have that
$m_{n,k-1}$, with $k=2,\ldots ,n$, is the number of permutations
of $S_{n}(1-23)$ ending with $k$ and $m_{n,n}$ is the number of
permutations of $S_{n}(1-23)$ ending with $1$. Then, if we move
the diagonal of $M$ such that it becomes the first column of the
matrix, we obtain a new matrix $A=(a_{i,j})_{i,j\geq 1}$ where
$a_{i,j}$ is the number of $(1-23)$-avoiding permutations of
length $i$ ending with $j$.

\scriptsize
$$ A =
\left (
\begin{array}{ccccccc}

                     1 & 0 & 0 & 0 & 0 & 0 & \vdots \\
                     1 & 1 & 0 & 0 & 0 & 0 & \vdots \cr
                     2 & 2 & 1 & 0 & 0 & 0 & \vdots \cr
                     5 & 5 & 3 & 2 & 0 & 0 & \vdots \cr
                     15 & 15 & 10 & 7 & 5 & 0 & \vdots \cr
                     52 & 52 & 37 & 27 & 20 & 15 & \vdots \cr
                     \cdots & \cdots & \cdots & \cdots & \cdots & \cdots & \ddots
\end{array}
\right)
$$

\normalsize The matrix $A$ is essentially the \emph{Bell
triangle}, which can be found in \cite{We} together with several
other references.

The above recursive properties of $M$ can be immediately
translated as follows:
\begin{enumerate}
    \item $a_{1,1}=1$ (the minimal permutation ends, trivially, with
    1);
    \item $a_{n,k}=0$ if $k>n$ (each permutation of length $n$ cannot end with a number greater than $n$
    itself);
    \item $a_{n,k}=\sum_{i=k}^{n-1}{a_{n-1,i}}+a_{n-1,1}$ if $2\leq k \leq n$ (the diagonal of $M$ has been moved in the first column of
    $A$);
    \item $a_{n,1}=a_{n,2}$ (since $a_{n,1}=m_{n,n}=m_{n,1}=a_{n,2}$).
\end{enumerate}
From 3 we obtain, for $k\geq 3$:
\begin{equation*}\label{eq:delta}
    a_{n,k}=a_{n,k-1}-a_{n-1,k-1},
\end{equation*}

If we denote by $\nabla$ the usual \emph{backward difference
operator}, since $a_{n,2}=B_{n-1}$, we get:

\begin{eqnarray*}
  a_{n,k} &=& \nabla a_{n,k-1} \\
    &=& \nabla^2 a_{n,k-1} \\
    &=& \cdots \\
    &=& \nabla^{k-2} a_{n,2}=\nabla^{k-2}B_{n-1}\quad \quad \mbox{(which holds also for $k=2$)}.
\end{eqnarray*}

Thus we find the following formulas concerning the distribution of
$1-23$-avoiding permutations according to their length and to the
value of their last entry:
$$
     | \{ \pi \in  S_{n}(1-23) : \pi_n = 1 \} | = B_{n-1} , \
    n \geq 1 ;
$$
$$
     | \{ \pi \in  S_{n}(1-23) : \pi_n = k \} | = \nabla^{k-2}
    (B_{n-1}) , \ 2 \leq k \leq n.
$$

\subsection{The other patterns of the class}\label{other}

The arguments employed for $S(1-23)$ can be easily modified for
the other patterns of the symmetry class of $1-23$, obtaining
similar results. The ECO construction, in these cases, has to be
adapted in order to obtain the same succession rule and the same
generating tree we got for $S(1-23)$. The matrices $M$ and $A$ are
defined as in the previous section.
\begin{enumerate}
    \item For the reverse pattern of $1-23$, i.e. $32-1$, we find
    that $a_{i,j}$ is the number of permutations $\pi$ of length $i$
        such that $\pi_1=j$, and so:
    \begin{itemize}
        \item $|\{\pi\in S_{n}(32-1):\pi_1=1\}|=B_{n-1},\
        n\geq 2$ ;
        \item $|\{\pi\in S_{n}(32-1):\pi_1=k\}|=\nabla^{k-2}
        (B_{n-1}),\ 2\leq k\leq n$.
    \end{itemize}
Note that in this case the ECO construction can be, in some way,
``reversed'', so that the active sites are not on the right of the
diagram of the permutation $\pi$ but on its left, i.e. before the
first entry of $\pi$.
    \item For the complement pattern $3-21$, we have that $a_{i,j}$
    is the number of permutations of length $i$ ending with $i+1-j$:
    \begin{itemize}
        \item $|\{\pi \in S_n(3-21):\pi_n=n\}|=B_{n-1}, n\geq 1$;
        \item $|\{\pi \in
        S_n(3-21):\pi_n=k\}|=\nabla^{n-k-1}(B_{n-1})$, \ $1\leq k \leq
        n-1$.
    \end{itemize}
        \item For the reverse-complement pattern $12-3$, $a_{i,j}$
        is the number of permutations $\pi$ of length
        $i$ such that $\pi_1=i+1-j$, and so:
    \begin{itemize}
        \item $|\{\pi\in S(12-3):\pi_1=n\}|=B_{n-1}$, \ $n\geq 1$;
        \item $|\{\pi\in
        S(12-3):\pi_1=k\}|=\nabla^{n-k-1}(B_{n-1})$, \ $1\leq k\leq
        n-1$.
    \end{itemize}
\end{enumerate}

\section{The symmetry class $\{3-12, 21-3, 1-32, 23-1\}$}

\subsection{ECO construction and generating tree of $S(3-12)$}

\begin{figure}
\begin{center}
\includegraphics[scale=0.4]{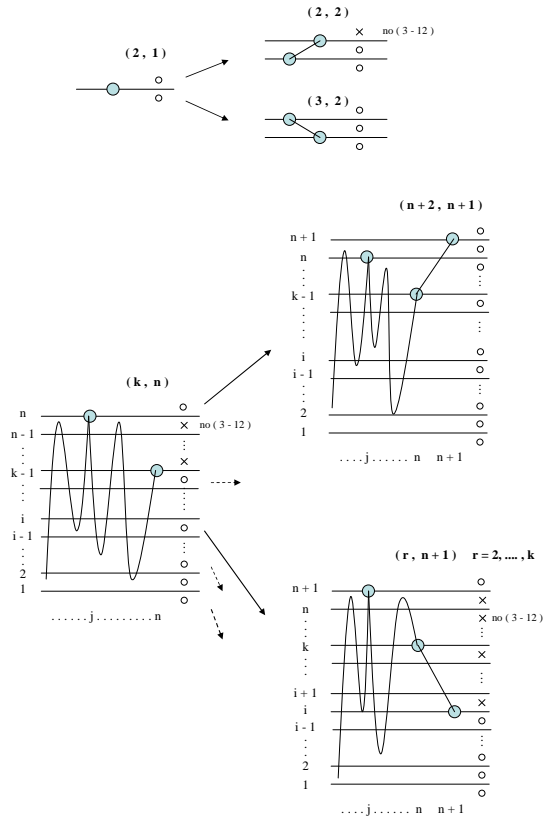}
\end{center}
\caption{ECO construction of $S(3-12)$}\label{fig:ECO3-12}
\end{figure}

Let $\pi\in S_n(3-12)$. If $\pi_n=k-1\neq n$, then $\pi$ generates
$k$ permutations $\pi^{(i)}\in S_{n+1}(3-12)$,
$i=1,2,\ldots,k-1,n+1$, by inserting a new node in region $i$. If
$\pi_n=n$, then $\pi$ generates $n+1$ permutations by inserting a
new node in any region. As it happened for the class $S(1-23)$,
note that the number of sons of $\pi$ is determined by the length
of $\pi$. It is easily seen that $\pi^{(n+1)}$ generates, in
turns, $n+2$ permutations, whereas $\pi^{(i)}$ ($i \neq n+1$)
produces $i+1$ permutations. This ECO construction is illustrated
in Figure \ref{fig:ECO3-12}. If each permutation of $S_n(3-12)$
with $k$ active sites is labelled $(k,n)$, then such a
construction can be encoded using the following succession rule:
$$
\left\{
\begin{array}{l}
  (2,1) \\
  (k,n) \rightsquigarrow (2,n+1) (3,n+1) \cdots (k,n+1)
(n+2,n+1) \quad .\\
\end{array}
\right.
$$

Since it is the same succession rule we got for $S(1-23)$, the
generating tree for $S(3-12)$ can be obtained in the same way.

\subsection{Distribution according to the length and the last value}

Defining the matrix $M=(m_{ij})_{i,j \geq1}$ as in Section
\ref{distr1-23}, it can be easily deduced that $m_{n,k}$ is the
number of permutations of $S_n(3-12)$ ending with $k$. Note that
in this case we do not need to move the diagonal of $M$ to obtain
the final matrix. Therefore, using again the backward difference
operator $\nabla$, the entries of $M$ have the form:
$$
m_{n,k}=\nabla^{k-1}(B_{n-1})
$$
whence:

$$|\{\pi\in S_{n}(3-12):\pi_n=n\}|=B_{n-1},\ n\geq 2\ ;$$
$$|\{\pi\in S_{n}(3-12):\pi_n=k\}|=\nabla^{k-1}(B_{n-1}),\ 1\leq
k\leq n-1\ .$$

\subsection{The other patterns of the class}
Proceeding as in Section \ref{other}, we get:
\begin{itemize}
    \item $|\{\pi\in S_{n}(21-3):\pi_1=n\}|=B_{n-1},\ n\geq 2$;
    \item $|\{\pi\in S_{n}(21-3):\pi_1=k\}|=\nabla^{k-1}(B_{n-1}),\ 1\leq k
    \leq n-1$;
\vspace{.1cm}
    \item $|\{\pi\in S_{n}(1-32):\pi_n=1\}|=B_{n-1},\ n\geq 2$;
    \item $|\{\pi\in S_{n}(1-32):\pi_n=k\}|=\nabla^{n-k}(B_{n-1}),\ 2\leq k\leq n$ ;
\vspace{.1cm}
    \item $|\{\pi\in S_{n}(23-1):\pi_1=1\}|=B_{n-1},\ n\geq 2$ ;
    \item $|\{\pi\in S_{n}(23-1):\pi_1=k\}|=\nabla^{n-k}(B_{n-1}),\ 2\leq k\leq n$.
\end{itemize}

\section{The symmetry class $\{2-13, 31-2, 2-31, 13-2 \}$}

The permutations of $S(2-13$) are enumerated by Catalan numbers
\cite{clae}. As far as the ECO construction of $S(2-13)$ is
concerned, we just note that, if $\pi\in S_n(2-13)$ is such that
$\pi_n=k$, then region $i$, for $i=1,2,\ldots,k+1$, is an active
site for $\pi$. The succession rule encoding this construction is:
$$
\left\{
\begin{array}{l}
  (2) \\
  (k) \rightsquigarrow (2)(3)\cdots(k+1)\\
\end{array}
\right.
$$

Defining the matrix $M$ as in the preceding sections, we obtain

\scriptsize
$$M=
\left(
\begin{array}{ccccccc}

                      1 & 0 & 0 & 0 & 0 & 0 & \vdots \\
                      1 & 1 & 0 & 0 & 0 & 0 & \vdots \cr
                      2 & 2 & 1 & 0 & 0 & 0 & \vdots \cr
                      5 & 5 & 3 & 1 & 0 & 0 & \vdots \cr
                      14 & 14 & 9 & 4 & 1 & 0 & \vdots \cr
                      42 & 42 & 28 & 14 & 5 & 1 & \vdots \cr
                     \cdots & \cdots & \cdots & \cdots & \cdots &\cdots & \ddots
\end{array}
\right)
$$
\normalsize which is the well-known Catalan Triangle whose entries
$m_{i,j}=\frac{j}{i}{2i-j-1\choose i-1}$ are the ballot numbers
and whose properties can be found, for example, in \cite{noo-zei}.

In the following, we present the results for all the patterns of
the class, which can be derived as in the previous sections (the
$m_{n,k}$'s are defined as before):
\begin{itemize}
    \item $|\{\pi\in S_{n}(2-13):\pi_n=k\}|=m_{n,k}=\frac{k}{n}{2n-k-1\choose n-1}$\ ;
    \item $|\{\pi\in S_{n}(31-2):\pi_1=k\}|=m_{n,k}=\frac{k}{n}{2n-k-1\choose n-1}$\ ;
    \item $|\{\pi\in S_{n}(2-31):\pi_n=k\}|=m_{n,n-k+1}=\frac{n-k+1}{n}{n+k-2\choose n-1}$\ ;
    \item $|\{\pi\in S_{n}(13-2):\pi_1=k\}|=m_{n,n-k+1}=\frac{n-k+1}{n}{n+k-2\choose n-1}$\ .
\end{itemize}

\section{Permutations avoiding a pair of generalized
patterns of type $(1,2)$ or $(2,1)$}

In  \cite{clae-man} Claesson and Mansour counted permutations
avoiding a pair of generalized patterns of type (1,2) or (2,1).
Similarly to what we have done in the previous sections, we can
study the distribution of the statistic ``first/last entry" on
permutations avoiding two or more generalized patterns. Here, we
consider only two special examples, the former being quite easy,
whereas the latter is surely more interesting. All the remaining
cases are left to the readers as open problems for future
research.

\subsection{An easy case}

We first deal with the permutations of $S(1-23, 1-32)$. This class
is enumerated by the number $I_n$ of involutions in $S_n$ (see
\cite{clae-man}). An ECO construction of this class can be encoded
by the following succession rule :
$$
\Omega: \left\{
\begin{array}{l}
  (2,1) \\
  (1,n) \rightsquigarrow (n+2,n+1) \\
  (n+1,n) \rightsquigarrow (1,n+1)^n (n+2,n+1) \\
\end{array}
\right.
$$
where the first element in the label is the number of active sites
of the permutation and the second one is its length. This can be
checked by representing permutations by means of the usual
path-like representation: indeed, if a permutation ends with $1$,
then an element can be inserted on its right in any region,
whereas if a permutation ends with $k\neq 1$, then the only
element which can be inserted must be placed in region $1$ on the
right. The reader is invited to complete the details, so to obtain
the construction described precisely by the succession rule
$\Omega$.

From the generating tree of $\Omega$, the matrix $M$ whose entry
$m_{i,j}$ is the number of vertices with label $j$ at level $i$
($i,j\geq 1$) can be constructed as in the preceding cases:

\scriptsize
$$
M= \left (
\begin{array}{ccccccccc}

                      0 & 1 & 0 & 0 & 0 &  0 &  0 &  0 & \vdots \\
                      1 & 0 & 1 & 0 & 0 &  0 &  0 &  0 & \vdots \cr
                      2 & 0 & 0 & 2 & 0 &  0 &  0 &  0 & \vdots \cr
                      6 & 0 & 0 & 0 & 4 &  0 &  0 &  0 & \vdots \cr
                     16 & 0 & 0 & 0 & 0 & 10 &  0 &  0 & \vdots \cr
                     50 & 0 & 0 & 0 & 0 &  0 & 26 &  0 & \vdots \cr
                    156 & 0 & 0 & 0 & 0 &  0 &  0 & 76 & \vdots \cr
                     \cdots & \cdots & \cdots & \cdots & \cdots & \cdots & \cdots & \cdots &\ddots
\end{array}
\right).
$$

\normalsize The entries can be immediately computed as follows:
\begin{itemize}
    \item $m_{1,1}=0$\ ,\quad $m_{1,2}=1$\ ;
    \item $m_{n,1}=(n-1)m_{n-1,n}$\ ,\quad $n\geq 2$\ ;
    \item $m_{n,n+1}=m_{n-1,1}+m_{n-1,n}$\ ,\quad $n\geq 2$\ ;
    \item $m_{i,j}=0$ in all the other cases.
\end{itemize}

From the ECO construction it easily appears that the first column
of $M$ counts the permutations $\pi$ of $S_n(1-23,1-32)$ such that
$\pi_{n-1}=1$ (or, which is the same, $\pi_n \neq 1$), whereas the
super-diagonal sequence $m_{n,n+1}$ ($n\geq 1$) shows the number
of $\pi$ ending with $1$. Since if $\pi\in S_n(1-23,1-32)$, then
$\pi_{n-1}=1$ or $\pi_n=1$, we deduce that the super-diagonal
satisfies $m_{n,n+1}=I_{n-1}$ ($n\geq 1$).

\subsection{A not so easy case}

\begin{figure}
\begin{center}
\includegraphics[scale=0.4]{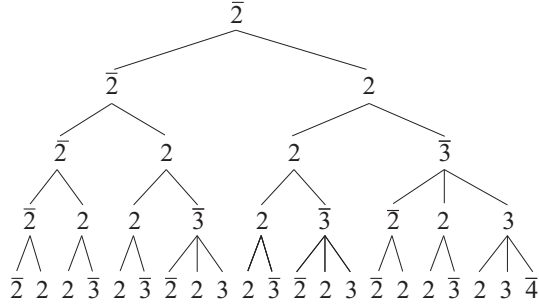}
\end{center}
\caption{the generating tree of
$S(1-23,21-3)$}\label{fig:tree123213}
\end{figure}

Our second example concerns the permutations of the class
$S(1-23,21-3)$, which also coincide with those of
$S(1-23,21-3,12-3)$(see \cite{ber-fer-pin}) and are enumerated by
Motzkin numbers. We will find the distribution of these
permutations according to their length and their last entry;
moreover, we will be able to derive the generating function of the
sequences enumerating the permutations of this class whose last
entry is $k$, for $k=1,2,\ldots$. We start by recalling the
coloured succession rule $\Phi$ encoding an ECO construction for
the above set of permutations (which can be found in
\cite{ber-fer-pin}):
$$
\Phi: \left\{
\begin{array}{l}
  (\bar{2}) \\
  (\bar{k}) \rightsquigarrow (\bar{2})(2)(3)\cdots(k)\\
  (k)\rightsquigarrow (2)(3)\cdots(k)(\overline{k+1}) \quad .\\
\end{array}
\right.
$$
In Figure \ref{fig:tree123213}, the first levels of the
corresponding generating tree are presented.

As in the preceding examples, we construct a matrix
$A=(a_{i,j})_{i,j\geq 1}$ recording in its entries the number of
labels at each level of the tree: namely, $a_{i,1}$ is the number
of coloured label $\bar{k}$, $k\geq 2$, at level $i$ of the tree
and $a_{i,j}$, $j\geq 2$, is the number of labels $j$ at level
$i$. The first lines of $A$ are:

\scriptsize
$$
A=
\left(%
\begin{array}{ccccccc}
  1 & 0 & 0 & 0 & 0 & 0 & \vdots \\
  1 & 1 & 0 & 0 & 0 & 0 & \vdots \\
  2 & 2 & 0 & 0 & 0 & 0 & \vdots \\
  4 & 4 & 1 & 0 & 0 & 0 & \vdots \\
  9 & 9 & 3 & 0 & 0 & 0 & \vdots \\
  21 & 21 & 8 & 1 & 0 & 0 & \vdots \\
  51 & 51 & 21 & 4 & 0 & 0 & \vdots \\
  127 & 127 & 55 & 13 & 1 & 0 & \vdots \\
  323 & 323 & 145 & 39 & 5 & 0 & \vdots\\
  835 & 835 & 385 & 113 & 19 & 1 & \vdots\\
  \cdots & \cdots & \cdots & \cdots & \cdots & \cdots & \ddots\\
\end{array}%
\right)\ .
$$

\normalsize As usual, we can find a recursive description of the
entries of $A$:
\begin{itemize}
    \item each label at level $i-1$ produces, among its sons, precisely one
    coloured label at level $i$, and so:
    $$a_{i,1}=\sum_{r\geq 1}a_{i-1,r}\quad ;$$
    \item each label $j\geq 2$ at level $i$ is generated either by a label $k\geq j$
    at level $i-1$ or by a coloured label $\bar{k}$, with $k\geq j$ at level $i-1$, which, in
    turn, is generated by the label $k-1$ at level $i-2$, then:
    \begin{equation}\label{ricorr}
        a_{i,j}=\sum_{k\geq j}a_{i-1,k}+\sum_{k\geq j-1}a_{i-2,k}\quad
        \mbox{for}\quad j\geq2\ ;
    \end{equation}

    \item it is easily seen that, in the above generating tree,
    the coloured label $\bar{k}$ first appears at the odd level
    $2k-3$, whereas the label $k$ first appears at the even
    level $2k-2$, whence: $$a_{i,j}=0\quad \mbox{for}\quad j\geq \lfloor i/2\rfloor +2\ .$$
\end{itemize}

The ECO construction of $S(1-23,21-3)$ shows that, if a
permutation has label $k$, then it ends with $k$, while if it has
a coloured label $\bar{k}$, then its last entry is $1$. Therefore,
the entry $a_{i,j}$ is the number of permutations with length $i$
and ending with the element $j$.

\bigskip

\noindent Our next aim is to find the generating function for the
sequences displayed in the columns of the matrix $A$, which are
the sequences enumerating the permutations of $S(1-23,21-3)$ with
last entry $j=1,2,\ldots$, according to their length. It is
convenient to change a little bit the notation: from now on, we
will index the lines of $A$ starting from $0$ instead of $1$.
First of all, we derive a simple recurrence for the entries of
$A$: using (\ref{ricorr}), simple calculations show that
\begin{equation}\label{ricorr2}
a_{n,k}=a_{n,k-1}-a_{n-1,k-1}-a_{n-2,k-2}\ , \quad \mbox{for}
\quad k\geq2 \ ,\ n\geq0\ .
\end{equation}
Let $C_k(x)$ be the generating function of the $k$-th column of
$A$:
$$
C_k(x)=\sum_{n\geq0} a_{n,k}x^n\ .
$$
Using (\ref{ricorr2}), we find the following recurrence relation
for $C_k(x)$:
\begin{equation}\label{CK}
C_{k+2}(x)=(1-x)C_{k+1}(x)-x^2C_k(x),\quad k\geq0\ .
\end{equation}
From the succession rule $\Phi$ (or from the ECO construction for
$S(1-23,21-3)$), it is easy to check that
$$
C_0(x)=M(x),\quad C_1(x)=M(x)-1\ ,
$$
where
$$
M(x)=\frac{1-x-\sqrt{1-2x-3x^2}}{2x^2}
$$
is the generating function of Motzkin numbers $\{M_n\}_{n\geq0}$.
In order to find a closed form for $C_k(x)$, we define a linear
operator $L$ on the vector space $\mathcal{V}$ of formal power
series of odd order. The set $(C_k(x))_{k\geq 1}$ is a basis of
$\mathcal{V}$, so $L$ can be defined as follows:
\begin{equation}\label{def_L}
L(C_k(x))=C_{k+1}(x)\quad\mbox{for}\quad k\geq1.
\end{equation}
From (\ref{CK}) it is:
$$
L^2(C_k(x))=(1-x)L(C_k(x))-x^2C_k(x)
$$
which is the same of
$$
(L^2-(1-x)L+x^2)C_k(x)=0\ .
$$
Therefore the operator $L^2-(1-x)L+x^2$ must vanish on
$\mathcal{V}$. Solving the equation $L^2-(1-x)L+x^2=0$, leads to
$$
L=\frac{1-x-\sqrt{1-2x-3x^2}}{2}=x^2M(x)\ .
$$
Now, from (\ref{def_L}), we obtain the desired closed form for
$C_k(x)$:
$$
C_k(x)=x^2M(x)C_{k-1}(x)=\cdots=x^{2(k-1)}M^{k-1}(x)(M(x)-1),\
k\geq1\ .
$$

\end{document}